\documentclass{birkjour}
\usepackage{amsfonts,amssymb,amsmath,mathtools,hyperref,dsfont}
\usepackage[english]{babel}
\usepackage{amssymb}

\DeclareMathAlphabet{\mathpzc}{OT1}{pzc}{m}{it}

\usepackage[normalem]{ulem}

\numberwithin{equation}{section}

\usepackage{graphicx,tikz}
\usepackage{mathrsfs}

\usepackage{color}

\newcommand{\ka}{\mathpzc k}
\newcommand{\el}{\mathpzc l}

\newcommand{\alg}{L^1(\R^+)}

\newtheorem{thm}{Theorem}[section]
\newtheorem{lemma}[thm]{Lemma}
\newtheorem{proposition}[thm]{Proposition}
\newtheorem{corollary}[thm]{Corollary}

\newcommand{\R}{{\mathbb R}}

\newcommand{\eps}{\epsilon}
\newcommand{\e}{\mathrm e}
\newcommand{\ud }{\, \mathrm d}

\newcommand{\dom}[1]{\mathcal D(#1)}
\newcommand{\grae}{\lim_{\eps \to 0+}}
\newcommand{\gra}{\lim_{n\to \infty}}
\newcommand{\lam}{\lambda}

\newcommand{\grat}{\lim_{t\to 0+}}
\newcommand{\grato}{\lim_{t\to \infty}}

\newcommand{\mc}{\mathcal}
\newcommand{\pol}{{\textstyle \frac 12}}
\newcommand{\mm}{\mathsf m}

\newcommand{\wt}{\widetilde}
\newcommand{\mquad}[1]{\qquad \text {#1} \qquad}
\newcommand{\mqquad}[1]{\quad \text {#1} \quad}
\newcommand{\ef}{{\sf \Sigma}}

\newcommand{\kap}{{\ka}}

\newcommand{\kamn}{{\el}}

\newcommand{\rod}[1]{\{#1 , t \ge 0 \}}

\newcommand{\alger}{L^1(\R)}
\newcommand{\kapa}{\ka^2}
\newcommand{\kapp}{\ka \times \ka}
\newcommand{\czn}{\int_0^\infty} 
\newcommand{\skap}{S_\ka}
\newcommand{\ma}{\mathfrak A}
\newcommand{\mg}{\mathfrak G_\eps}
\newcommand{\Psik}{L^1(\skap)}

\title[Approximation of Walsh's spider process]
{A kinetic model approximation of Walsh's spider process on the infinite star-like graph\\ \begin{center} \small{To the memory of Jan Kisyński (1933--2022)}\end{center}}

\begin{document}

\author{Adam Bobrowski}
\email{a.bobrowski@pollub.pl}
\address{Lublin University of Technology, Nadbystrzycka 38A, 20-618 Lublin, Poland.}

\author{El\.zbieta Ratajczyk}
\email{e.ratajczyk@pollub.pl}
\address{Lublin University of Technology, Nadbystrzycka 38A, 20-618 Lublin, Poland.}


\date{\today {\bf File: {\jobname}.tex.}} 


\begin{abstract}We consider processes of deterministic motions on $\ka$ copies  of the star-like graph $\skap\coloneqq K_{1,\ka}$ with $\ka$ edges which are perturbed by two stochastic mechanisms: one caused by interfaces located at the graphs' centers, the other describing jumps between different copies of the same edge. We prove, extending the main result of \cite{abtk}, that diffusing scaling of these processes leads in the limit to the Walsh's spider process on $\skap$.  \end{abstract}

\thanks{Version of \today}


\maketitle

\vspace{-1cm}
\section{Introduction}

\subsection{From telegraph equation to a L\'evy process on a non-commutative group and Wiener process}
\begin{figure}
\begin{tikzpicture}[scale=0.7]
\draw [thick,->] (-5,0) -- (5,0);
\node [above] at  (4.8,0.1) {\tiny{{$\mathbb R\times\{1\}$}}}; 
\draw [thick,<-] (-5,-2) -- (5,-2);
\node [above] at  (4.8,-1.9) {\tiny{{$\mathbb R\times\{-1\}$}}}; 
\draw [magenta,<->,dashed] (3,0)--(3,-2); 
\draw [magenta,<->,dashed] (-3,0)--(-3,-2); 
\end{tikzpicture}
\caption{\footnotesize{$\mathbb G$ is a Lie group with natural motion along the curves $t \mapsto (vt,1) \circ g$ starting from $g\in \mathbb G$, that is, the motion to the right with speed $v$ when started at the upper copy $\R\times \{+1\}\subset \mathbb G$, and the motion to the left with speed $v$ on the lower copy $\R \times \{-1\}\subset \mathbb G$.  In the process $\Lambda$ this motion is perturbed by jumps from one copy to the other at the epochs of the Poisson process.}}
\label{rys1}
\end{figure}
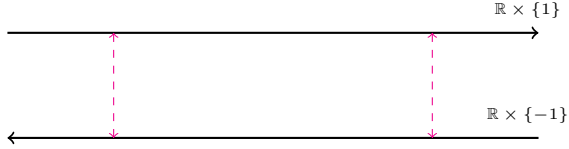

Let $a$ and $v$ be two positive constants. S. Goldstein \cite{goldsteins} was apparently the first to notice that the telegraph equation 
\begin{equation} \label{intro:1} \partial_{tt} u (t,x) + 2a \partial_t u (t,x) = v^2 \partial_{xx}u (t,x), \qquad t\ge 0, x \in \R, \end{equation}
a hyperbolic PDE by nature, exhibits properties that are usually considered to be attributes of  parabolic PDEs of special type, that is, of Kolmogorov equations for Markov processes. In particular, the Cauchy problem for \eqref{intro:1} with initial condition 
\begin{equation}\label{intro:2} u (0,x) = u_0(x), \quad \partial_t u (0,x) = 0, \qquad x \in \R, \end{equation} 
 is well posed and its solution is nonnegative whenever $u_0$ is. Later, M. Kac \cite{kac} expressed this solution in terms of the stochastic process  underlying  \eqref{intro:1} as follows: 
\begin{equation}\label{intro:3} u (t,x) = \pol \mathbb E \, [ u_0(t+v\xi_a (t)) + u_0(t-v\xi_a(t))],  \qquad t \ge 0, x \in \R,  \end{equation} 
where $\mathbb E$ stands for expected value, 
\begin{equation} \xi_a (t) \coloneqq \int_0^t (-1)^{N_a(s)} \ud s, \qquad t \ge 0, \label{intro:4} \end{equation}
and $N_a(t), t \ge 0$ is the Poisson process with $\mathbb E N_a(t) = at, t \ge 0$. M. Kac's seminal paper, in turn, opened the way for  the development of the theory of \emph{random evolutions} of Griego and Hersh  
\cite{gh1,gh2}, see also  \cite{ethier} Chapter 12 and \cite{pinskyrandom}. J. Kisyński \cite{kkac} chose a slightly different direction and has shown that the possibility of expressing solutions of \eqref{intro:1}--\eqref{intro:2} in the form \eqref{intro:3} hinges on the fact that  (see Figure \ref{rys1})
\begin{equation}\Lambda (t) \coloneqq (v \xi_a(t), (-1)^{N_a(t)} ), \qquad t \ge 0 \label{intro:5} \end{equation}
is a L\'evy process with values in the locally compact, non-commutative group \[ \mathbb G \coloneqq \R \times \{-1,1\}\] with multiplication defined by $(x,k)\circ (y,\ell) = (x\ell +y, k\ell)$; for the general theory of such processes see \cite{heyer}. Markovian nature of $\Lambda $ is also crucially used in the exposition of the telegraph equation in \cite{ethier}.

Assuming $a=v^2$ and letting $a \to \infty$, we see that, at least formally, in the limit, \eqref{intro:1} becomes the diffusion equation
\begin{equation} \label{intro:6} \partial_t u (t,x) = \pol \partial_{xx}u (t,x), \qquad t\ge 0, x \in \R. \end{equation}
This heuristic reasoning may be made precise: in can be  proved that solutions of \eqref{intro:1}  converge to those of \eqref{intro:6} (see e.g. \cite{banbob,deg,knigazcup,ethier} and references given there). From the perspective of processes, this limit theorem can be interpreted as follows. By letting $a\to \infty$ we make the jumps from one part of $\mathbb G$ to the other (see Figure \ref{rys1} again) so frequent that, in the limit, two points: $(x,1)$ and $(x,-1)$, are lumped into one for all $x\in \R$, and thus the limit state-space is not $\mathbb G$ but $\R$. Moreover, as suggested by \eqref{intro:6}, the limit process is a Brownian motion. See e.g. \cite{jansen} and the already cited \cite{pinskyrandom} for more on this subject.

\newcommand{\cer}{\mathfrak  C[-\infty,\infty]}

\subsection{Perturbation at an interface leads to skew Brownian motion} 

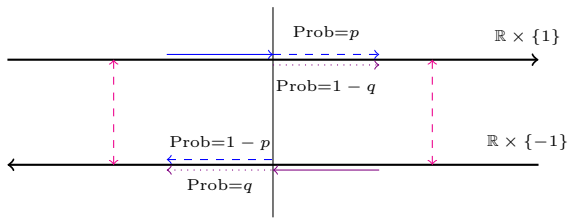
\begin{figure}
\begin{tikzpicture}[scale=0.7]
\draw [thick,->] (-5,0) -- (5,0);
\node [above] at  (4.8,0.1) {\tiny{{$\mathbb R\times\{1\}$}}}; 
\draw [thick,<-] (-5,-2) -- (5,-2);
\node [above] at  (4.8,-1.9) {\tiny{{$\mathbb R\times\{-1\}$}}}; 
\draw [thin] (0,-3)--(0,1); 
\draw [blue,->] (-2,0.1) -- (0,0.1);
\draw [blue,dashed,->] (0,0.1) -- (2,0.1);
\draw [blue,dashed,<-] (-2,-1.9) -- (0,-1.9);
\node [above] at  (1,0.2) {\tiny{Prob=$p$}}; 
\node [above] at  (-1,-1.9) {\tiny{Prob=$1-p$}}; 
\draw [violet,dotted,<-] (-2,-2.1) -- (0,-2.1);
\draw [violet,<-] (0,-2.1) -- (2,-2.1);
\draw [violet,dotted,->] (0,-0.1) -- (2,-0.1);
\node [below] at  (1,-0.2) {\tiny{Prob=$1-q$}}; 
\node [below] at  (-1,-2.1) {\tiny{Prob=$q$}};  
\draw [magenta,<->,dashed] (3,0)--(3,-2); 
\draw [magenta,<->,dashed] (-3,0)--(-3,-2); 
\end{tikzpicture}
\caption{\footnotesize{A perturbation of the process $\Lambda$ of \eqref{intro:5} and Figure \ref{rys1}.  Particles moving to the left (on the lower line) and to the right (on the upper line),  may be reflected at the interface with probabilities depending on wether  they approach the interface from the left or from the right.}}
\label{rys2}
\end{figure} 

Much more recently, in the paper \cite{abtk}, inspired partly by the kinetic model of a motion of a phonon with an interface,  studied in
\cite{tomekkinetic,tomekkinetic3,tomekkinetic2,tomekkinetic1}, and the telegraph
process   with elastic boundary at the origin \cite{wlosi,wlosi1}, it has been discovered that  by introducing an additional perturbing mechanism at an interface one can alter the limit process discussed above: the regular Brownian motion becomes \emph{skew} Brownian motion.  

The latter,  introduced in \cite{ito,walsh},  is a natural generalization of the standard one-dimensional Brownian motion: it behaves like a Brownian motion except that the sign of each excursion from $0$ is chosen using an independent Bernoulli random variable --- see \cite{lejayskew} for much more information on the process. The skew Brownian motion turns out to be an honest Feller process on $\R$, and as such can be described by means of 
its generator, that is, a Laplace operator, say, $A$, in the space $C_0(\R)$  of continuous functions on $\R$ that vanish at infinity. The domain of $A$ is composed of functions satisfying the following three properties: 
\begin{itemize}
\item[(a)] $f$ is twice continuously differentiable in both
  $(-\infty,0]$ and $[0,\infty)$, separately, with left-hand and
  right-hand derivatives at $x=0$, respectively, 
\item[(b)] $\lim_{x\to \infty} f'' (x)=\lim_{x\to -\infty} f''(x) =0$,  
 \item[(c)] for certain positive $\alpha$ and $\beta$ boundary conditions 
\begin{equation}\label{intro:7} f''(0+)=f''(0-) \mquad{ and } \beta f'(0-)=\alpha f'(0+),\end{equation} 
  hold; note that the first of them implies that, although $f'(0)$ need not exist, it is meaningful to speak of $f''(0)$.  
\end{itemize}
Furthermore, we
   define \( A f \coloneqq \pol f'' .\)
 Parameters $\alpha$ and $\beta$ have the following interpretation: $\frac{\alpha}{\alpha + \beta}$ is the probability that the sign of excursion is chosen to be positive, and $\frac{\beta}{\alpha + \beta}$ is the probability that the sign is negative.

The interface alluded to above, which changes the limit standard Brownian motion to the skew Brownian motion, is located 
 at the  points {$(0,\pm1)\in \mathbb G$} and  works as follows (see Figure \ref{rys2}). A~particle obeying the rules of the process $\Lambda$ of \eqref{intro:5} and approaching  the interface from
the left, thus {moving} on the upper copy, filters through the
interface with probability $p$ and continues its motion to the right on $\R \times
\{1\}$. With probability $1-p$, however, the particle is
reflected and starts moving to the left (from $(0,-1)$) on the lower
copy. Analogously, when approaching the interface from the right (on the lower copy), the particle filters through the interface with probability $q$ or is reflected and continues its motion on the upper copy with probability $1-q$.  One of the main results of \cite{abtk} says that, provided that $a=v^2$, the so-perturbed process $\Lambda$ converges, as $a\to \infty$, to the skew Brownian motion with parameters $\alpha=p$ and $\beta=q$.

\subsection{The goal of the paper}\label{sec:tgoo}
\begin{figure}
\begin{tikzpicture}
      \node[circle,fill=blue,inner sep=0pt,minimum size=2pt] at (360:0mm) (center) {};
    \foreach \n in {1,...,8}{
        \node [circle,fill=blue,inner sep=0pt,minimum size=0.1pt] at ({\n*360/8}:2cm) (n\n) {};
        \draw [blue,dashed](center)--(n\n);}
     \node[circle,fill=blue,inner sep=0pt,minimum size=2pt] at (360:0mm) (center) {};
    \foreach \n in {1,...,8}{
        \node [circle,fill=blue,inner sep=0pt,minimum size=0.1pt] at ({\n*360/8}:1.5cm) (n\n) {};
        \draw [blue,thick](center)--(n\n);}   
                 \node [above] at (4,-0.3) {$\displaystyle\sum_{i=1}^\ka \alpha_i f_i'(0)=0$}; 
        \node [below] at (0,1.5) {\includegraphics[scale=0.165]{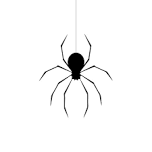}};\end{tikzpicture}
\caption{\footnotesize{The infinite star-like graph $\skap$  with $\kap=8$ edges. 
Walsh's process on $\skap\coloneqq K_{1,\ka}$ is a Feller process whose behavior at the graph's center is characterized by the boundary condition visible above; outside of the center the process behaves like a standard Brownian motion.}}\label{slg} 
\end{figure}
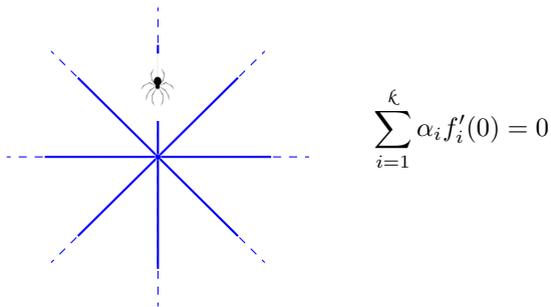

Among many generalizations of the skew Brownian motion, one that seems to have attracted particular attention is the Walsh's spider process on the infinite star-like graph $\skap\coloneqq K_{1,\ka}$ (see e.g. see \cite{kostrykin2012,lejayskew,manyor,yor97}) with $\kap$ edges --- see Figure \ref{slg}. This process is characterized by positive parameters $\alpha_1,\dots, \alpha_{\kap}$ such that $\sum_{i=1}^\kap \alpha_i =1 $, playing the role of probabilities. Roughly, when at the graph's center, Walsh's spider chooses the $i$th edge with probability $\alpha_i$ to continue its motion there; outside of the center it follows the rules of  a standard one-dimensional Brownian motion (see Section \ref{sec:srt} for more information).    

The goal of our paper is to find an approximation of Walsh's spider process by means of processes analogous to those considered in \cite{abtk}. Certainly, the skew Brownian motion is a particular case of Walsh's spider, corresponding to $\ka =2$. Since to obtain in the limit the skew Brownian motion, a process on $S_2$, one needs to consider approximating processes with values on two copies of $S_2$,  it seems reasonable
to look for approximations of the Walsh's spider process in $\skap$ among the processes with values on $\kap$ copies of $\skap$.

The approximating processes we construct are mixtures of two simpler ones, say, $X$ and $Y$. To describe the first of these we imagine (see Figure \ref{slgmotion}) a particle which, when on the $i$th edge of an $i$th copy of $\skap$, moves deterministically towards the graph's center with speed $\ka -1 $. The center is an interface which introduces randomness to the motion.
This means that after reaching the interface the particle 
\begin{enumerate} 
\item either continues its motion on the same copy of the graph, choosing the $j$th edge with probability $p_{i,j}$; it then moves away from the center with speed $1$, 
\item or, with probability $p_{i,i}\not =0$, jumps to another copy of the graph $\skap$; 
 conditional on such a jump the probability of choosing a $j$th copy is $r_{i,j}$; then the particle moves away from the center on the $i$th edge of the $j$th copy of the graph.  
\end{enumerate}
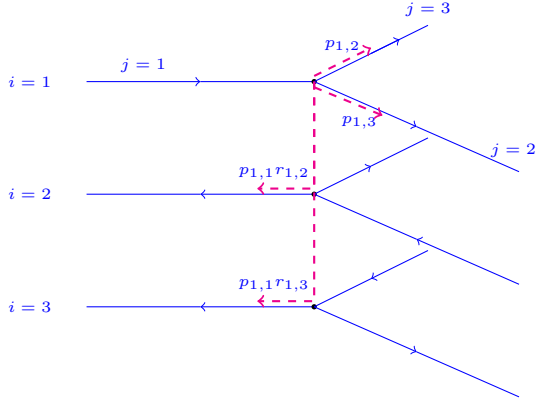
\begin{figure}
\begin{tikzpicture}[color=blue,scale=0.75]
\begin{scope}
\node at (-5,0) {\tiny{$i=1$}};  \node [above] at (-3,0) {\tiny{$j=1$}};   
\node [above] at (3.5,-1.5) {\tiny{$j=2$}};   \node [above] at (2,1) {\tiny{$j=3$}};
\node [above] at (0.5,0.35) {\tiny{$p_{1,2}$}};  \node [above] at (0.8,-1) {\tiny{$p_{1,3}$}}; 
\draw [->] (-4,0)--(-2,0);\draw (-2,0)--(0,0);
\node[circle,fill=black,inner sep=0pt,minimum size=2pt] at (360:0mm) (center) {};
\draw [->] (0,0)--(1.8,-0.8);\draw [dashed,->,color=magenta,thick] (0,-0.1)--(1.2,-0.6);
\draw (1.8,-0.8)--(3.6,-1.6);
\draw [->] (0,0)--(1.5,0.75); \draw [dashed,->,color=magenta,thick] (0,0.1)--(1,0.6);
\draw [-] (1,0.5)--(2,1);
\end{scope}
\begin{scope}[shift={(0,-2)}]
\node at (-5,0) {\tiny{$i=2$}}; 
\draw [-] (-4,0)--(-2,0);\draw [<-] (-2,0)--(0,0);
\node[circle,fill=black,inner sep=0pt,minimum size=2pt] at (360:0mm) (center) {};
\draw [-] (0,0)--(1.8,-0.8);
\draw [<-] (1.8,-0.8)--(3.6,-1.6);
\draw [->] (0,0)--(1,0.5);
\draw [-] (1,0.5)--(2,1);
\end{scope}
\begin{scope}[shift={(0,-4)}]
\node at (-5,0) {\tiny{$i=3$}}; 
\draw [-] (-4,0)--(-2,0);\draw [<-] (-2,0)--(0,0);
\node[circle,fill=black,inner sep=0pt,minimum size=2pt] at (360:0mm) (center) {};
\draw [->] (0,0)--(1.8,-0.8);
\draw (1.8,-0.8)--(3.6,-1.6);
\draw [-] (0,0)--(1,0.5);
\draw [<-] (1,0.5)--(2,1);
\end{scope}
\draw [dashed,color=magenta,thick] (0,0)--(0,-1.9)[->]--(-1,-1.9);
\draw [dashed,color=magenta,thick] (0,0)--(0,-3.9)[->]--(-1,-3.9);
\node [above] at (-0.7,-1.9) {\tiny{${p_{1,1}}r_{1,2}$}}; 
\node [above] at (-0.7,-3.9) {\tiny{${p_{1,1}r_{1,3}}$}}; 
\end{tikzpicture}
\caption{\footnotesize{Scattering at the graph's center: A particle moves towards the graph center at the first edge of the first copy of the graph; there, with probability $p_{1,2}$ continues its motion on the second edge of the same graph, or, with probability $p_{1,3}$, on the third copy of the same graph. With probability $p_{1,1}$, however, it is `reflected' and starts moving away from the center on the first edge of  the second or third copy of the graph; conditional on reflection the probabilities of choosing the second and the third copies are $r_{1,2}$ and $r_{1,3}$, respectively.}} 
\label{slgmotion}
\end{figure} 

In particular, the particle moves towards the center (with speed $\ka -1$) only on the $i$th edge of an $i$th copy of the $\skap$ graph, $i =1,\dots, \kap$; on all the remaining edges it moves away from the center with speed $1$. The probabilities $p_{i,j}$ and $r_{i,j}$ (by convention, $r_{i,i}=0, i=1,\dots, \kap$) form two $\kapp $ transition probability matrices
\begin{equation}\label{intro:8} P = (p_{i,j})_{i,j=1,\dots, \kap} \quad \text{and} \quad R = (r_{i,j})_{i,j=1,\dots, \kap} \end{equation}
that is, matrices of non-negative numbers in which elements in each row add up to $1$. 

The second component, $Y$,  of the approximating processes is a random scattering mechanism playing the role of 
jumps between the lower and upper copies of $\R$ in the process $\Lambda$ of \eqref{intro:5}, depicted in Figures \ref{rys1} and \ref{rys2} by dashed lines. Namely, a particle moving on a copy of the $j$th edge of $\skap$ will at random times, as governed by a time-continuous Markov chain's intensity matrix $Q^j$, jump to the $j$th edge of another copy of $\skap$ without changing its distance from the origin. 

In the main theorem of the paper, Thm \ref{thm:2}, we show that there is a relatively large class of intensity matrices $Q^j, j=1,\dots, \ka$ with the following property: given parameters $\alpha_i, i=1,\dots, \kap$, there is a 
family of interrelated probability matrices $P$ and $R$ such that the two processes described above, when  combined and appropriately scaled, converge to the Walsh's spider process on $\skap$. 

Our proof is based on the theory of convergence of semigroups of operators, as expounded e.g. in \cite{kniga,knigazcup,ethier}. Semigroups that describe the approximating processes are presented in Section \ref{semigroups}, whereas the semigroup that describes the Walsh's process is presented in Section \ref{sec:srt}. 
Section \ref{sec:aow} discusses the details of the main approximation theorem and its assumptions. The proof of the theorem is contained in Section \ref{pota}; all the necessary lemmas are gathered in Section \ref{fkl}.

\vspace{-0.5cm}
\section{Semigroups that describe approximating processes}\label{semigroups}

\subsection{The space} 
Let $\alg$ and $\alger$ be the spaces of (classes of) absolutely integrable functions on $\R^+\coloneqq [0,\infty)$ and $\R$, respectively. 
We start by considering the Cartesian product 
\[ \Phi \coloneqq [\alg]^{\kapa}\] 
of $\kapa $ copies of $\alg$. A member of $\Phi$ can thus be seen as a $\kapp$ matrix of elements $\phi (\cdot, i, j)$ of $\alg$, where $i,j$ belong to the set 
\[ \mc K \coloneqq \{1,\dots, \ka \}.\] 
For the norm in $\Phi$ we choose
\( \|\phi\| \coloneqq \sum_{i,j\in \mc K} \|\phi (\cdot, i, j)\|_{\alg}.\)
This space is isometrically isomorphic to the space of integrable functions on $\ka$ copies of the star-like graph $\skap= K_{1,\kap}$ with $\kap$ edges (see Figure \ref{slg}) in which all (infinitely long) edges emanating from the graph's center are identified with the half-line $\R^+$, equipped with the one dimensional Lebesgue measure.  
In other words, each $\phi \in \Phi$ can be identified with a single function on $\kap$ copies of the graph  $\skap$; then $i$ is the number of the copy of  $\skap$  and $j$ is the number of the edge in $\skap$.

\newcommand{\stargraph}[2]{\begin{tikzpicture}
      \node[circle,fill=black,inner sep=0pt,minimum size=2pt] at (360:0mm) (center) {};
    \foreach \n in {1,...,#1}{
        \node at ({\n*360/#1}:#2cm) (n\n) {}; 
        \draw (center)--(n\n);} 
\end{tikzpicture}}


\subsection{Markov semigroup for component $X$}\label{sec:dod}


 With the help of transition matrices $P$ and $R$ of \eqref{intro:8} we define a strongly continuous semigroup of operators in $\Phi$ as follows:
\begin{equation}\label{dod:1} T(t) \phi (x, i, i) = \phi (x+\kamn t, i, i) \qquad i\in \mc K \end{equation}
and, for $i\not =j$,  
\begin{equation}\label{dod:2} T(t)\phi (x, i, j) =\begin{cases} \phi (x-t, i,j), & x \ge t, \\
\kamn p_{i,j} \phi (\kamn (t-x),i,i) +\el p_{j,j}r_{j,i} \phi (\kamn (t-x), j,j), & x < t ;\end{cases} \end{equation}
here and in what follows, to simplify and shorten formulae, 
\begin{equation} \el \coloneqq \ka - 1.\label{def_l}\end{equation}
A straightforward calculation establishes the semigroup property \[ T(t)T(s)=T(t+s), \qquad s,t\ge 0,\] and the fact that each $T(t)$ is a Markov operator. The latter statement means that $\sum_{i,j\in \mc K} \czn T(t) \phi (x,i,j) \ud x = \sum_{i,j\in \mc K} \czn \phi (x,i,j) \ud x $ and  $T(t)\phi (\cdot, i, j) \ge 0$, provided that we have $\phi (\cdot, i, j)\ge 0$ for $i,j\in \mc K$. It follows that $\|T(t)\|=1$, that is, that $\rod{T(t)}$ is a semigroup of contractions.

This semigroup describes the component $X$ of the approximating process, as introduced in Section \ref{sec:tgoo}. By  this we mean that if $\phi$ is an initial distribution of $X$, then $T(t)\phi$ is its distribution at time $t\ge 0.$ More precisely, $T(t)\phi (\cdot, i,j)$ is the density of the probability that $X$ at time $t$ is at the $j$th edge of the $i$th copy of the graph $\skap$.


We claim that the semigroup $\rod{T(t)}$ is strongly continuous and that its generator, say, $\mc A$, is characterized as follows.

\begin{proposition}\label{prop1} 
The domain $\dom{\mc A}$ of $\mc A$ is composed of $\phi \in \Phi$ such that: \begin{itemize}
\item [(a) ] Each $\phi (\cdot, i,j), i,j\in \mc K$ is absolutely continuous with absolutely integrable derivative, that is, there are $C_{i,j}\in \R, i,j\in \mc K$ and a $\varphi \in \Phi$ such that 
\begin{equation}\label{sctg:1}\phi (x,i,j) = C_{i,j} + \int_0^x \varphi (y, i,j) \ud y, \qquad x\ge 0, i,j\in \mc K. \end{equation}
\item [(b) ] The following transmission conditions are satisfied: 
\begin{equation}\label{sctg:2} \phi (0,i,j) = \kamn p_{i,j} \phi (0,i,i) + \el p_{j,j}r_{j,i} \phi (0,j,j),  \qquad i\not = j; i,j\in \mc K,\end{equation}
that is, $C_{i,j}$s  of \eqref{sctg:1} are interrelated as follows 
\( C_{i,j} = \kamn p_{i,j} C_{i,i} + \el p_{j,j}r_{j,i}C_{j,j}. \)
\end{itemize}
Moreover, for such $\phi$ we have
\begin{align} \mc A \phi (\cdot, i, i)& =\ \kamn \phi' (\cdot, i,i) = \, \kamn \varphi (\cdot, i,i), \qquad i \in \mc K, \nonumber  \\
 \mc A \phi (\cdot, i, j) &= - \phi' (\cdot, i,j)= -\varphi (\cdot, i,j),  \qquad i\not = j; i,j\in \mc K. \label{sctg:dod} \end{align}
\end{proposition}


\begin{proof}
To begin, we consider the Cartesian product  
\[ \wt \Phi \coloneqq [\alg ]^\kap \times [\alger]^{\ka (\ka -1)}.\]
As in the case of $\Phi$, we think of a member $\wt \phi$ of $\wt \Phi$ as an $\kapp$ matrix. However, now the diagonal elements belong to $\alg$ whereas the off-diagonal elements belong to $\alger$:  
\[ \wt \phi (\cdot, i,i) \in \alg, \quad i\in \mc K \mquad{and} \wt \phi (\cdot, i,j) \in \alger, \quad i\not =j;i,j\in \mc K.\]
Moreover, we define the semigroup $\rod{\wt T(t)}$ in $\wt \Phi$ by 
\begin{equation}\label{sctg:3} \wt T(t)  \wt \phi (x, i, i) = \wt \phi (x+\kamn t, i, i) \qquad i\in \mc K \end{equation}
and, for $i\not =j$,  
\begin{equation}\label{sctg:4} \wt T(t)\phi (x, i, j) = \wt \phi (x-t, i,j). \end{equation}
This is just the Cartesian product semigroup built of left translations in $\alg$ and right translations in $\alger$. Hence, it is obviously strongly continuous and its generator $\wt {\mc A}$ is characterized as follows (see e.g. \cite[pp. 66--67]{engel}). A $\wt \phi$ belongs to $\dom{\wt{\mc A}}$ if there is a $\wt \varphi \in \wt \Phi$ and real constants $\wt C_{i,j}$ such that  
\begin{equation}\label{sctg:5}\wt \phi (x,i,i) = \wt C_{i,i} + \int_0^x \wt \varphi (y, i,j) \ud y, \qquad x\ge 0,\end{equation}
and 
\begin{equation}\label{sctg:6}\phi (x,i,j) = \wt C_{i,j} + \int_0^x \wt \varphi (y, i,j) \ud y, \qquad x\in \R, i\not =j; i,j\in \mc K. \end{equation}
For such $\wt \phi$, 
\[ \wt{\mc A} \wt \phi (\cdot, i, i) = \kamn \wt \varphi (\cdot, i,i), \, i \in \mc K \mqquad {and}
\wt{\mc A}\wt \phi (\cdot, i, j) =  -\wt \varphi (\cdot, i,j),  \, i\not = j; i,j\in \mc K. \] 

Next, we consider the subspace $\wt \Phi_0 $ of $\wt \Phi $, composed of $\wt \phi \in \wt \Phi$ such that 
\begin{equation}\label{sctg:7} \wt \phi (-x, i,j) = \kamn p_{i,j} \wt \phi (\kamn x,i,i) + \el p_{j,j}r_{j,i} \wt \phi (\kamn x, j,j), \end{equation} for $ x \ge 0, i\not = j; i,j\in \mc K$.
It is a key observation, checked by a straightforward calculation, that $\wt \Phi_0$ is invariant under $\rod{\wt T(t)}$. The family $\rod{[\wt T(t)]_{|\wt \Phi_0}}$ is thus a strongly continuous semigroup in $\wt \Phi_0$ (termed the subspace semigroup, see \cite{engel}) and its generator $\wt {\mc A}_0$ is the restriction of $\wt {\mc A}$ to the domain $\dom{\wt {\mc A}_0}\coloneqq \dom {\wt {\mc A}} \cap \wt \Phi_0$. 

Finally, we observe that $\Phi$ is isomorphic to $\wt \Phi_0$; this is just to say that an entire matrix $\wt \phi \in \wt 
\Phi_0$ is determined by its diagonal entries plus the restrictions of its off-diagonal entries to $\R^+$. More formally, the operator $\mc E\colon\Phi \to \wt \Phi$ given by 
\[ \mc E \phi (x , i, j) = \phi (x, i, j) \qquad x\ge 0, i,j \in \mc K\]
and                                
\[ \mc E \phi (-x , i, j) = \kamn p_{i,j} \wt \phi (\kamn x,i,i) + \el p_{j,j}r_{j,i} \wt \phi (\kamn x, j,j), \qquad  x \ge 0, i\not = j; i,j\in \mc K\]
is linear, injective, and maps $\Phi$ onto $\wt \Phi_0$. Moreover,
\begin{align*}
\|\mc  E\phi\|_{\wt \Phi} &= \|\phi\|_{\Phi} + \sum_{i\in \mc K} \sum_{j\not =i} [ p_{i,j} \czn |\phi (x,i,i) |\ud x + p_{j,j}r_{j,i} \czn |\phi (x,j,j)| \ud x ]  \\
& =  \|\phi\|_{\Phi} +  \sum_{i\in \mc K} \czn |\phi (x,i,i) |\ud x \le 2 \|\phi\|_{\Phi}, 
 \end{align*}
so that $\mc E$ is bounded with norm $2$ (the upper bound is obtained whenever $\phi$ vanishes outside the main diagonal) and $\mc E$ has a left and right inverse $\mc R$, where $\mc R \colon\wt \Phi_0 \to \Phi$ is the restriction operator 
\[  \mc R \wt \phi (x,i,j) = \wt \phi (x,i,j), \qquad x\ge 0, i,j\in \mc K.\]

The discussed objects are related to $\rod{T(t)}$ by the following formula
\begin{equation}\label{sctg:8} T(t)\phi = \mc R \wt T(t)\mc E \phi = \mc R [\wt T(t)]_{\wt \Phi_0} \mc E\phi , \qquad \phi \in \Phi, t \ge 0;\end{equation} 
this means that $\rod{T(t)}$ in $\Phi$ is isomorphic to the subspace semigroup $\rod{[\wt T(t)]_{\wt \Phi_0}}$ in $\wt \Phi_0$. 
It follows that $\rod{T(t)}$ is a strongly continuous semigroup, for so is $\rod{[\wt T(t)]_{\wt \Phi_0}}$. Moreover, a $\phi\in \Phi$ belongs to $\dom{\mc A}$ iff $\mc E\phi$ belongs to $\dom{\wt{\mc A}_0} = \dom{\wt {\mc A}}\cap \wt \Phi_0$. It is now easy to check that this is the case iff conditions (a) and (b) of the definition of $\dom{\mc A}$ are satisfied. Formula  \eqref{sctg:8} implies also
\begin{align*} \mc A \phi (x,i,i) & = \mc R\wt{\mc A}_0 \mc E\phi (x,i,i) = \kamn \phi'(x,i,i), \ \, \qquad x\ge 0, i \in \mc K,\\
 \mc A \phi (x,i,j) & =\mc  R\wt{\mc A}_0 \mc E\phi (x,i,j) = - \phi'(x,i,j), \qquad x\ge 0, i\not =j; i,j \in \mc K. \end{align*}
This completes the proof. \end{proof}

\subsection{Markov semigroup for component $Y$} 

Component $Y$ is a process of jumps between the same edges of different copies of $\skap$: while on the $i$th copy of the $j$th edge of $\skap$ a particle may jump to the $k$th copy of the $j$th edge of $\skap$ (without changing its distance from the origin) as  in a continuous time Markov chain governed by an intensity matrix 
\begin{equation}\label{msf:1} Q^j=\big(q^j_{i,k}\big)_{i,k \in \mc K}.\end{equation}
In other words, the scattering mechanism is governed by a family of Markov chains which perturb the argument $i$ of $\phi \in \Phi$
while keeping   $j$ the same.

The related semigroup of Markov operators is generated by the following bounded linear operator in $\Phi$:
\[  \mc Q\phi (\cdot, i, j) = \sum_{k\in \mc K} q^j_{k,i}\phi (\cdot, k,j), \qquad i,j\in \mc K.  \]
The exponential function of $\mc Q$, for $t\ge 0$,  is given by 
\[ \e^{t\mc Q} \phi (\cdot, i,j) = \sum_{k \in \mc K} p^j_{k,i} (t) \phi (\cdot, k,j), \qquad i,j\in \mc K,   \] 
where $(p_{k,i}^j(t))_{k,i\in \mc K}$ is the transition probability matrix for the Markov chain with intensity matrix \eqref{msf:1}.

\subsection{Generators of the diffusing scaling processes}

It is the subject of this paper to study the limit of a diffusing scaling of the `mixture' of the processes described above. In other words, we want to find the limit, as $\eps \to 0+$, of the semigroups in $\Phi$ generated by 
\begin{equation} \mc G_\eps \coloneqq  \eps^{-1} \mc A + \eps^{-2} \mc Q;\label{gotd:0}\end{equation}
the fact that each $\mc G_\eps $ is a Markov semigroup generator can be proven as in Section 2.2 in \cite{abtk}. 
We will show that, under certain conditions on $Q^j$, $P$ and $R$ (see Section \ref{sec:aow}), the semigroups generated by $\mc G_\eps$ converge to the Markov semigroup describing Walsh's spider process on $\skap$ --- see Theorem \ref{thm:2} for a precise statement.

\section{Semigroups related to the Walsh's spider process on $\skap$}\label{sec:srt}

\subsection{Walsh's spider process as a Feller process}
Let, as in Section \ref{sec:tgoo}, $\alpha_j, j\in \mc K$ be positive numbers adding up to $1$. It will be convenient to assume, without loss of generality, that the sequence $\alpha= (\alpha_j)_{j \in \mc K}$ is ordered, that is, that $\alpha_1\le \alpha_2 \le \dots \le \alpha_{\kap}$. As shown in \cite{barlow}, the Walsh's spider process  with characteristic $\alpha$ is a Feller process on $\skap$. Furthermore, the related semigroup $\rod{\mc T_\alpha (t)}$
of operators in $C_0(\skap)$, the space of continuous functions on $\skap$ that vanish at infinity, can be given rather explicitly by  means of the semigroups describing the minimal Brownian motion on   $\skap$ and the reflecting Brownian motion on $[0,\infty)$ --- see \eqref{wsp:1} below. 

The minimal Brownian motion on $\skap$, while on one of the edges, away from the graph center $O$, behaves like a standard one-dimensional Brownian motion. However, at the first moment it touches $O$, it is killed and removed from the state-space. Strictly speaking, thus, its state-space is not $\skap$ but $\skap^0 \coloneqq \skap \setminus \{O\}$. Now, any member $f \in C_0(\skap^0)$ can be identified with the sequence $(f_j)_{j\in \mc K}$ of elements of the space $C_0(0,\infty)$ of continuous functions on the positive half-line that vanish at both  $0$ and  $\infty$.  Moreover,  the minimal Brownian motion semigroup $\rod{\mc T_{\text{min}} (t)}$ on $\skap$ can be identified with the Cartesian product of $\ka$ copies of the familiar minimal Brownian motion semigroup  $\rod{T_{\text{min}} (t)}$ in $C_0(0,\infty)$.  It follows that the domain of the generator, say, $G_{\text{min}}$, of  $\rod{\mc T_{\text{min}} (t)}$, is composed of $(f_j)_{j\in \mc K}$ such that all $f_i\in C_0(0,\infty)$ are twice continuously differentiable with $f_i''\in C_0(0,\infty)$, and   $G_{\text{min}}(f_j)_{j\in \mc K}=\pol (f_j'')_{j\in \mc K}$. 

We recall also that the reflecting Brownian motion on $[0,\infty)$ starting at an $x\in [0,\infty)$ is defined as $|x + w(t)|, t \ge 0$, where $w(t), t \ge 0$ is a standard one-dimensional Brownian motion on $\R$ starting at $0$. The related semigroup in the space $C_0[0,\infty)$ of continuous functions on $[0,\infty)$ that vanish at infinity is given by 
\( T_{\text{ref}} (t) f (x) = \mathbb E \, f (|x + w(t)|), t,x \ge 0, f \in C_0[0,\infty).\)
The domain of the generator $G_{\text{ref}}$ of $\rod{ T_{\text{ref}} (t) }$ is composed of $f$ that are twice continuously differentiable with $f'' \in C_0[0,\infty)$, and satisfy $f'(0)=0$; for such $f$ we have  $ G_{\text{ref}}f=\frac 12 f''$.

To express $\rod{\mc T_\alpha(t)}$, a semigroup in $C_0(\skap)$,  by means of  the semigroups $\rod{\mc T_{\text{min}} (t) }$  and  $\rod{ T_{\text{ref}} (t) }$ we note finally that any $f\in C_0(\skap)$ can be identified with a sequence $(f_j)_{j\in \mc K}$ of elements of $C_0[0,\infty)$ such that $f_i (0)=f_j(0), i,j\in \mc K$. With this identification in mind, for $f= (f_j)_{j\in \mc K} $, we have (see \cite{barlow}, eq. (2.2)) 
\begin{equation}\label{wsp:1}  \mc T_\alpha (t) f = \mc T_{\text{min}} (t) (f - \mc I \overline f ) + \mc I T_{\text{ref}} (t) \overline f, \qquad f\in C_0(\skap), t \ge 0, \end{equation}
where $\overline f \coloneqq \sum_{j=1}^\ka \alpha_j f_j \in C_0[0,\infty)$ and $\mc I\colon C_0[0,\infty) \to C_0(\skap)$ assigns the constant sequence $(g)_{j\in \mc K}$ to a $g\in C_0[0,\infty)$.

The generator $G_\alpha$ of  $\rod{\mc T_\alpha(t)}$ is characterized as follows: its domain is composed of $(f_j)_{j\in \mc K}$ such that each $f_j$ is twice continuously differentiable with $f''_j \in C_0[0,\infty)$, and 
\begin{equation}\label{azda} \overline f' (0)= \sum_{j\in \mc K} \alpha_j f'_j(0)=0 \mquad{ and } f''_i(0)=f_j''(0), \quad i,j\in \mc K;\end{equation} for such $(f_j)_{j\in \mc K}$, \( G_\alpha (f_j)_{j\in \mc K}=\pol f'' \coloneqq  \pol (f_j'')_{j\in \mc K}.\) Indeed, on one hand, for $f$ described above, $f- \mc I \overline f $ belongs to the domain of $G_{\text{min}} $ and $\overline f$ belongs to the domain of $G_{\text{ref}}$. Therefore, 
\[ \grat t^{-1} (\mc T_\alpha (t) f - f) = G_{\text{min}} (f - \mc I \overline f) 
+ \mc I G_{\text{ref}}\overline f=\pol (f_j'')_{j\in \mc K} \in C_0(\skap), \]
showing that the generator of $\rod{\mc T_\alpha(t)}$ extends $G_\alpha$. On the other hand, calculating as in \cite{kostrykin1} and \cite{konkusSIMA} one can check that given a $\lam >0$ and a $g\in C_0(\skap)$ there is precisely one $f \in \dom{G_\alpha}$ such that $\lam f - G_\alpha f = g$. A standard argument shows thus that the searched for generator cannot be a proper extension of $G_\alpha$. A different derivation of the boundary conditions \eqref{azda} can be found in \cite{kostrykin2012}.
\vspace{-0.3cm}
\subsection{The `dual' Markov semigroup}
The component $X$ of the approximating processes does not posses the Feller property. As a result, the approximating semigroups generated by operators $\mc G_\eps$ of \eqref{gotd:0} are not defined in a space of continuous functions, but in the $L^1$ type space $\Phi$ of Section \ref{semigroups}. We cannot thus expect that in the limit the semigroup  $\rod{\mc T_\alpha(t)}$ of \eqref{wsp:1} will be obtained. Rather, we should expect the limit semigroup to be in a sense dual to  $\rod{\mc T_\alpha(t)}$. 
Here are the details. 

First of all, we equip $\skap$ with the measure, say, $\mm$, which at each of the edges coincides with the usual Lebesgue measure. By the Riesz representation theorem, the space $L^1(\skap)$ of functions on $\skap$ that are integrable with respect to $\mm$ can be seen as a subset of the dual $[C_0(\skap)]^*$. On the other hand, $L^1(\skap)$ can be identified with  the Cartesian product of $\kap$ copies of $\alg$: 
\vspace{-0.1cm}
\begin{align}\label{wsp:2}L^1(\skap) \overset{id}\coloneqq [\alg]^\kap,\end{align}
Thus, as in the case of $\Phi$, a $\psi \in \Psik$ has a dual status: it can either be seen as a vector $\left ( \psi (\cdot, j)\right )_{j\in \mc K}$ of elements of $\alg$,  or a single function on the $S_\kap$ graph;  $j$ is thought of as the index of the graph's edge. The space  $\Psik$ is equipped  with the usual norm
\( \|\psi\|=\int_{\skap} |\psi| \ud \mm  = \sum_{j\in \mc K} \|\psi (\cdot,j)\|_{\alg}  \).

In $\Psik$ we define an operator $A_\alpha$ as follows. Its domain is composed of $\psi \in \Psik$ such that 
\begin{itemize} 
\item [(a) ] $\psi (\cdot, j) \in W^{2,1} (\R^+), j\in \mc K$, that is, for each $j\in \mc K$, $\psi (\cdot,j)$ is twice differentiable with $\psi''(\cdot, j)$ in $\alg$, and  
\item [(b) ] there are constants $C$ and $D_j,j\in \mc K$ such that 
\[ \psi (x, j) = \alpha_j C + D_j x + \int_0^x (x-y) \psi''(y,j) \ud y, \qquad x\ge 0, j\in \mc K, \]
whereas $\sum_{j\in \mc K} D_j=0$; in other words, $\psi \in \dom{A_\alpha}$ satisfy the following transmission conditions 
\begin{equation}\label{wsp:3} \sum_{j\in \mc K} \psi'(0,j) =0 \quad \text{ and } \quad \alpha_j^{-1} \psi (0,j) = \alpha_i^{-1} \psi (0,i), \quad i,j\in \mc K. \end{equation}
 \end{itemize}
For such $\psi$ we let \( A_\alpha \psi = \pol \psi''\coloneqq \pol (\psi''(\cdot, j))_{j\in \mc K}.\)

The following proposition reveals a connection between operators $A_\alpha$ and $G_\alpha$. It says that the dual  
$\rod{\mc T_\alpha^* (t)}$ to $\rod{\mc T_\alpha (t)}$ leaves the subspace $L^1(\skap) \subset [C_0(\skap)]^*$ invariant, forms a strongly continuous semigroup of operators there, and as restricted to this subspace is generated by $A_\alpha$. 
It means in particular that Walsh's spider process, besides having Fellerian nature, has the following property: if its initial distribution is absolutely continuous with respect to $\mm$ then so is its distribution at all $t\ge 0$. If $\psi$ is its initial density, then $\e^{tA_\alpha}\psi $ is its density at time $t\ge 0$.

\begin{proposition} $A_\alpha$ is the generator of a semigroup of Markov operators in $L^1(\skap)$. Moreover, for $f\in \dom{G_\alpha} $ and $\psi \in \dom{A_\alpha}$, 
\[ \int_{\skap}  f A_\alpha \psi  \ud \mm = \int_{\skap} (G_\alpha f)\psi \ud \mm. \] 
\end{proposition}

\begin{proof}Since $A_\alpha$ is clearly densely defined and a short calculation using \eqref{wsp:3} establishes that $\int_{\skap} A_\alpha \psi \ud \mm =0$ for all $\psi \in \dom{A_\alpha}$, to prove the first sentence,  by  \cite{rudnickityran} Thm. 4.4., we need to check that for any $\psi \in L^1(\skap)$ and $\lam > 0$ there is a unique $\psi_0 \in \dom{A_\alpha}$ such that $\lam \psi_0 - A_\alpha \psi_0 = \psi$; moreover, $\psi_0\ge 0$ whenever $\psi\ge 0$.

\newcommand{\slam}{\sqrt {2\lam}}
Such a $\psi_0$ has to be of the form 
\begin{align} 
\psi_0(x,j) &= C_j \e^{\slam x} + D_j \e^{-\slam x}
 - \sqrt{\frac 2\lam}\int_0^x \sinh \slam (x-y) \psi(y,j) \ud y \nonumber \\ &=
\frac 1{\slam}\int_0^\infty \e^{-\slam |x-y|} \psi(y,j) \ud y + D_j \e^{-\slam x},\  x \ge 0, j\in \mc K, \label{wsp:4} \end{align}
where $C_j \coloneqq \frac 1{\slam}\int_0^\infty \e^{-\slam y} \psi(y,j) \ud y$ (otherwise, $\psi_0 (\cdot, j)$ is not integrable) and 
$D_j$ are to be determined. Since $\psi_0'(0,j) = \slam (C_j -D_j)$, the boundary conditions \eqref{wsp:3} are satisfied iff 
\begin{equation} \sum_{j\in \mc K} C_j =  \sum_{j\in \mc K} D_j \mqquad {and} \alpha_j^{-1} (C_j+D_j) = \alpha_i^{-1} (C_i +D_i), \quad i,j \in \mc K.\label{wsp:5} \end{equation}
This system, in turn, has the unique solution 
\begin{equation} D_j = 2\alpha_j \sum_{k\in \mc K}C_k  - C_j, 
\qquad j \in \mc K.\label{wsp:6} \end{equation}
For this choice of constants, $C_j +D_j \ge 0, j \in \mc K$  whenever $\psi\ge 0$. Hence, 
\begin{align*} \psi_0 (x,j) &\ge 
\frac 1{\slam}\int_0^\infty \e^{-\slam |x-y|} \psi(y,j) \ud y - C_i \e^{-\slam x}
\\
& \ge  \frac 1{\slam} \int_0^\infty \left [ \e^{-\slam |x-y|} - \e^{-\slam (x+y)} \right ]\psi (y,j) \ud y \ge 0,\end{align*}
as long as $\psi \ge 0$, completing the proof of the first part. The rest is established by a straightforward calculation. \end{proof}
 

\newcommand{\ulam}{\textstyle{\frac 1\ka}}

\section{Approximation of Walsh's spider process}\label{sec:aow}

\subsection{Choice of $Q$}\label{sec:coq}

To mimic the properties of the model discussed in \cite{abtk} in our more general situation, we 
assume that each $Q^j$ is symmetric, and 
\begin{equation} q^j_{j,j}=-\ka + 1 \mquad { and } q^j_{j,i}=1, \quad i\neq j, j\in \mc K.\label{sco:1}\end{equation} 
It follows that each $Q^j$, $j \in \mc K$ is  an irreducible intensity matrix, having the vector \( \frac 1{\ka} (1,1,\dots, 1)\) as invariant distribution. Thus, by  \cite{norris} Thm 3.6.2, we see that  $\grato p_{k,i}^j (t)= \frac 1\kap$ for $i,j,k\in \mc K$; this in turn renders 
\begin{equation}\label{sco:2} \grato \e^{t\mc Q} \phi (\cdot, i, j) = \ulam \sum_{k\in \mc K}  \phi (\cdot, k,j), \qquad i,j\in \mc K.   \end{equation}
We note that, since $Q^j$ is symmetric, $\e^{tQ^j}$ is doubly stochastic: in each row its elements add up to $1$ and so do its elements in each column.  

\subsection{Choice of $P$ and $R$}\label{sec:cop}

It is intuitively clear that different choices of matrices $P$ and $R$ of \eqref{intro:8} lead to different limits for the semigroups \eqref{gotd:0}, or no reasonable limit at all. To say the least, given $\alpha$ as above we should not expect that all choices of $P$ and $R$ will lead to the Walsh's spider process with this particular parameter. It turns out (see Lemma \ref{lem:tsos1} further down) that if we want our approximation scheme to work, we should restrict ourselves to matrices $P$ and $R$ related by the following constrains: 
\begin{equation}\label{tsos:1}\alpha_j (1-(\ka-1)p_{j,j}r_{j,i}) = (\ka -1)\alpha_i p_{i,j}, \qquad i\not =j;i,j\in \mc K. \end{equation} 
They play a somewhat similar role to the \emph{detailed balance conditions} (see e.g., \cite[p.~48]{norris} or \cite[p. 322]{spectralgap}), and in particular imply that $\alpha$ 
is an invariant measure for $P$.

The following example shows that the family of such pairs of matrices is non-empty. Given $\ka \ge 3$ and $\delta \in [\frac 1{\ka - 1}, 1]$, we define
\[ r_{1,2} =\delta, r_{1,j} = {\textstyle \frac{1-\delta}{\ka -2}}, \quad j =3,\dots, \ka \mqquad{ and }r_{i,j} = {\textstyle \frac 1{\ka -1}}, \quad  i\not =1, j \not = i;\]
(as always $r_{i,i}=0, i\in \mc K$). This form of $R$ forces the off-diagonal terms of the related $P$ to be
\[ p_{2,1}= \tfrac{\alpha_1 (1-(\ka -1)\delta p_{1,1})}{\alpha_2 (\ka -1)},  \qquad p_{i,1}= \tfrac{\alpha_1 (1-\frac{\ka -1}{\ka -2}(1-\delta) p_{1,1})}{\alpha_i (\ka -1)}, \quad i=3,\dots, \ka,\]
and $p_{i,j} = \frac {\alpha_j (1-p_{j,j})}{\alpha_i (\ka -1)},i\not =j, j=2,\dots,\ka$. Hence, the question of existence of a $P$ that is related to $R$ via \eqref{tsos:1} reduces to that of whether non-negative $p_{i,i}, i \in \mc K$ can be chosen in such a way that  the above formulae define a transition matrix. We claim that $P$ is such a matrix if 
\[ p_{1,1} =  \tfrac{\ka (\ka-2)(1-\gamma)}{(\ka -1 -\delta)(\ka-1)}, p_{2,2}=1 - \tfrac {\alpha_1}{\alpha_2}  
\tfrac{(\ka -2)(\delta+1)\gamma - (\ka-1)\delta +1}{\ka -1 -\delta}, p_{i,i} = 1- \tfrac{\alpha_1}{\alpha_i} \gamma, i\ge 3\]
for a $\gamma \in [\gamma_0, 1]$ where $\gamma_0 = \gamma_0 (\delta) = \frac{(\ka-1)((\ka-1)\delta -1)}{\ka (\ka -2)\delta}$
(in fact, it can be proved that the above formulae give the general solution to the problem of existence of $P$ for the particular $R$ under consideration). To show this  we check first that $0\le p_{1,1}\le \frac 1{(\ka -1)\delta} \le 1 \le \frac {\ka -2}{(\ka-1)(1-\delta)}$ (where, by convention $\frac {k-2}0\coloneqq \infty$) and $p_{i,i} \in [0,1]$ for $i=2,\dots, \ka$; it follows that all $p_{i,j}$s are non-negative. Moreover, {a longer calculation confirms that elements in each row of $P$ add up to $1$.}



\newcommand{\staroc}{
\begin{proposition} A probability matrix $P$ satisfies \eqref{tsos:1} iff there is a $\gamma \in [0,1]$ such that 
\begin{equation}\label{tsos:2} p_{i,i} = 1 - {\textstyle \frac {\alpha_1}{\alpha_i}}\gamma,\quad p_{i,j} ={\textstyle \frac {\alpha_1}{(\ka -1)\alpha_i}}\gamma, \qquad j\not =i, i,j\in \mc K. \end{equation}
\end{proposition}

\begin{proof} Let $\gamma_j \coloneqq 1-p_{j,j}, j \in \mc K$. Conditions \eqref{tsos:1} are satisfied iff $p_{i,j} = \frac {\alpha_j}{(\ka -1)\alpha_i} \gamma_j, i \not =j; i,j\in \mc K$. Also, since $\sum_{j\in \mc K} p_{i,j}$ is to equal $1$, we must have $\alpha_i \gamma_i = \frac 1{\ka -1} \sum_{j\not =i} \alpha_j \gamma_j, i \in \mc K$. This, however, means that $\alpha_i \gamma_i = \alpha_j \gamma_j, i,j\in \mc K$, that is, that \eqref{tsos:2} holds (with $\gamma\coloneqq \gamma_1$).  Finally, all quantities in \eqref{tsos:2} lie in $[0,1]$ iff 
so does $\gamma$, and then $\sum_{i\in \mc K} p_{i,j} =1$. Hence, it remains to check that $p_{i,j}$ defined in \eqref{tsos:2} satisfy \eqref{tsos:1}. \end{proof}
In particular, in a general $2\times 2$ transition probability matrix with non-zero entries
\( \begin{pmatrix} 1-p & p\\ q & 1-q\end{pmatrix},\) where $p,q\in (0,1)$,  we can always arrange (by possibly exchanging its rows) that $q\ge p$, to see that it is of the form \eqref{tsos:2} with, for example, $\alpha_1 =p$ and $\alpha_2 = q$. Thus, in the case $\ka =2$ conditions \eqref{tsos:2} are automatically satisfied for certain $\alpha_1,\alpha_2$; for this reason there was no need to study them in \cite{abtk}.  

In terms of  $\gamma_i \coloneqq 1-p_{i,i}, i \in \mc K$, \eqref{tsos:2} can be written as 
\begin{equation}\label{tsos:2a} p_{i,i} = 1 - \gamma_i,\quad p_{i,j} ={\textstyle \frac {\gamma_i}{\ka -1}} \qquad j\not =i, i,j\in \mc K; \end{equation} 
in particular, $p_{i,j}, j\not =i$ depend only on $i$.}

A couple of remarks are worth making here. First of all, for $\ka =2$, there is only one possible $R$, that is, $R= {\tiny \begin{pmatrix}0 & 1\\ 1 & 0\end{pmatrix}}$. Furthermore, in a general $2\times 2$ transition probability matrix with non-zero entries
\( {\tiny \begin{pmatrix} 1-p & p\\ q & 1-q\end{pmatrix}},\) where $p,q\in (0,1)$,  we can always arrange (by possibly exchanging its rows) that $q\le p$, to see that \eqref{tsos:1} holds with $\alpha_1 =\frac q{p+q}$ and $\alpha_2 = \frac p{p+q}$. Thus, in the case $\ka =2$ conditions \eqref{tsos:1} are automatically satisfied; for this reason there was no need to study them in \cite{abtk}. Secondly, for $\delta =\frac 1{\ka -1}$, $P$ has a particularly simple form 
\( p_{i,i} = 1 - {\textstyle \frac {\alpha_1}{\alpha_i}}\gamma, p_{i,j} ={\textstyle \frac {\alpha_1}{(\ka -1)\alpha_i}}\gamma,\) \(j\not =i, i,j\in \mc K, \)
where $\gamma \in [0,1]$. Thirdly, the class of pairs of $P$ and $R$ related via \eqref{tsos:1} is apparently much larger than that discussed above: in fact, for a number of randomly chosen matrices $R$, Maple was able to find a corresponding $P$. However, the  problem of determining all $P$ and $R$ related via \eqref{tsos:1} exceeds the scope of this paper.  
  
Returning to the two-parameter family of matrices $P$ and $R$ we note that, except for the case $\gamma =\gamma_0$, all off-diagonal entries in $P$ are positive. It follows that the discrete-time Markov chain with transition probability matrix $P$ is irreducible and aperiodic, and $\alpha \in \R^{\kap}$  is its invariant measure. Hence, see e.g. \cite[p. 41, Thm. 1.8.3]{norris} or \cite[p. 310, Thm. 1.2.1]{spectralgap}, \begin{equation}\label{ergo}\gra P^n = \Pi \end{equation} 
where $\Pi = (\pi_{i,j})_{i,j\in \mc K} \in \R^{\ka^2}$ is defined by $  \pi_{i,j} = \alpha_j, i,j\in \mc K$.
Our main theorem holds under assumption that both \eqref{tsos:1} and \eqref{ergo} are satisfied.

\subsection{The approximation theorem} \label{sec:tat}

The space 
$\Psik$ of \eqref{wsp:2} is isometrically isomorphic to the subspace $\Phi_0$ of
$\Phi$ made of functions $\varphi $ such that $\varphi (\cdot, i,j)$ does not depend on $i$.  The isomorphism we have in mind is $J\colon\Psik \to \Phi_0$ given
by 
\begin{equation}
J\psi (\cdot ,i,j) = \ka^{-1} \psi (\cdot ,j), \quad  i,j\in \mc K, \psi \in \Psik,\label{tat:1}
\end{equation}
{with $J^{-1}\varphi (\cdot ,j) = \ka \varphi (\cdot, 1,j), \varphi \in \Phi_0$.} 
 It follows that the operators
\[ S(t) \coloneqq J \e^{t A_\alpha} J^{-1},\qquad  t \ge 0 \]
form a strongly continuous semigroup of operators in $\Phi_0$. Its
generator is 
\begin{equation}
\label{tat:2}\widetilde A_\alpha \coloneqq J A_\alpha
J^{-1},
\end{equation}
 with the domain equal to  
$\dom{\widetilde A_\alpha} =
J\big(\dom{A_\alpha}\big)$.
That is, a $\varphi \in \Phi_0$ belongs to
$\dom{\widetilde A_\alpha} $ iff $J^{-1} \varphi $ is in $\dom{A_\alpha}$
and then  $\widetilde A_\alpha \varphi = JA_\alpha J^{-1} \varphi  =\pol \varphi''$, see e.g. \cite[Section 7.4.22]{kniga}.

\newcommand{\dab}{A_\alpha}
\newcommand{\wdab}{\widetilde{A_\alpha}}

\begin{thm}\label{thm:2}  Assume that the matrices $Q^j$, $P$ and $R$ are as described in Sections \ref{sec:coq} and \ref{sec:cop}. Let $\mc P\colon\Phi \to \Phi$ be the projection on $\Phi_0$ defined by \(\mc P \phi =  \left (\tfrac 1\kap{\textstyle \sum_{i\in \mc K}} \phi (\cdot, i, j) \right )_{i,j\in \mc K}, \phi \in \Phi.\) 
Then, for $c\coloneqq 2 \frac {\ka -1}\ka$,
\[ 
{\grae} \e^{t \mc G_\eps } \phi = \e^{c t \widetilde {\dab}} \mc P \phi, \qquad t >0, \phi \in \Phi, \]
{strongly in the norm of $\Phi$}, and the limit is uniform in
$t$ on compact subsets of $(0,\infty)$. For $\phi \in \Phi_0$, the
limit holds also for $t=0$ and is uniform in $t$ on compact subsets of $[0,\infty)$. 
\end{thm}
The explain the meaning of this theorem, let us think of a probability density $\phi\in \Phi$ (i.e., a non-negative 
function of  norm $1$), interpreted as the density of the initial distribution of the process generated by $\mc G_\eps$. Our theorem says that, for any $t>0$, as $\eps \to 0$,   $\e^{t\mc G_\eps } \phi $   loses its dependence on $i$, and in the limit 
can be identified (via $J$) with a member of $\Psik$ which is the density of the Walsh's spider process at $t$, provided that the initial density of the Walsh's process is $J^{-1}\mc P\phi \in \Psik$.

The theorem will be proved in Section \ref{pota}; Section \ref{fkl} gathers all the necessary lemmas. We note that, besides Kurtz's singular perturbation theorem  (see \cite{ethier,kurtzper,kurtzapp} or
  \cite[Thm. 42.1]{knigazcup}), our argument involves the ideas of \cite{greiner}
(see also Theorem 3.1 in \cite{marta+g}),
\cite[pp. 230--232]{banbob2}  and \cite[Lemma 2.3]{banasiak}

\vspace{-0.4cm}
\section{Four key lemmas}\label{fkl}

Our first lemma characterizes the kernel of $\lam - \mg $ for $\lam >0$ and $\eps >0$, where $\mg$ is an extension of 
$\mc G_\eps$ of \eqref{gotd:0} defined as follows. First, we enlarge $\mc A$ to the operator $\ma$ with domain 
$\dom{\ma}\supset \dom{\mc A}$ composed of $\phi\in \Phi$ of the form \eqref{sctg:1}, where constants $C_{i,j}$ need not satisfy transmission conditions \eqref{sctg:2}. Moreover, for such $\phi$ we agree that, as in \eqref{sctg:dod}, 
\begin{align*} \ma \phi (\cdot, i, i)& =\ \el \phi' (\cdot, i,i) = \, \el \varphi (\cdot, i,i), \qquad i \in \mc K, \\
 \ma \phi (\cdot, i, j) &= - \phi' (\cdot, i,j)= -\varphi (\cdot, i,j),  \qquad i\not = j; i,j\in \mc K, \end{align*}
where, to recall $\el =\ka -1$.  Then, we define  
\[ \mg \coloneqq \eps^{-1} \ma + \eps^{-2} \mc Q , \qquad \eps >0.\]

Notations of the lemma involve $\el \times \el $ matrices 
\(\widetilde{Q}^j\coloneqq (\tilde{q}^j_{i,k})_{ i,k\in \mc K\setminus\{j\}}, j \in \mc K,\) where
$  \tilde{q}^j_{i,i}={q}^j_{i,i}+1$ and $  \tilde{q}^j_{i,k} ={q}^j_{i,k}$ for $i\neq k$. In other words, $\widetilde{Q}^j$ is obtained by removing the $j$th row and $j$th column of the matrix $Q^j+{I}$, where $I$ is the $\kapp$ identity matrix. Because of assumption \eqref{sco:1}, $\widetilde Q^j$ is a (symmetric) intensity matrix. We write 
\[ \e^{t \widetilde{Q}^j}\coloneqq (\tilde{p}^j_{i,k} (t) )_{ i,k\in \mc K\setminus\{j\}}, \qquad t\geq 0,   \] 
to denote the related doubly stochastic matrix of transition probabilities.

Finally, given $\lam >0$ and $\eps >0$, we define $\mu =\mu (\eps, \lam)$ and $\nu= \nu (\eps, \lam)$  by  
\begin{equation}\label{gotd:1} \mu  \coloneqq \tfrac {\eps \lam (\ka -2) + \sqrt{(\lam \eps \ka)^2 + 4\ka(\ka-1)\lam}}{2(\ka-1)}, \text{ } \nu \coloneqq  \lam \eps +\eps^{-1}. \end{equation} 
Direct calculations
 verify that \( \tfrac{1}{\eps(\nu-\mu)}=\eps \mu +\tfrac{\eps^2\lam}{\ka-1}+1\), 
implying, in particular, that  $\nu>\mu$.

\begin{lemma}\label{lem:gotd1} Let $\eps,\lam >0$ be given. A $\phi\in \Phi $ belongs to the kernel of $\lam - \mg$ iff there are real constants $E_{i,j}, i,j\in \mc K$ such that, for $ x \ge 0, i,j \in \mc K,$
\begin{align} 
\phi (x,j,j) & = E_{j,j} \e^{-\mu x}, \nonumber \\ 
\phi (x,i,j) & = \e^{-\nu x} \sum_{k\neq j} E_{k,j}   \tilde{p}^j_{k,i}(x/\eps )+ \tfrac{E_{j,j}}{\eps (\nu - \mu)} (\e^{-\mu x}-\e^{-\nu x}),\quad i \neq j, \label{gotd:3}\end{align}
for $\mu$ and $\nu$ introduced above, and $E_{i,j}$s satisfy \( \sum_{i\not =j} E_{i,j} =\tfrac{(\ka-1) E_{j,j}}{\eps(\nu-\mu)}, j\in \mc K.\)  \end{lemma}

\begin{proof}A $\phi\in \dom{\mg}$ belongs to the kernel of $\lam - \mg$ iff for all $j\in \mc K$,
\begin{align}
(\lam \eps^2+ \el ) \phi (\cdot, j, j) - \eps \el \phi '(\cdot, j, j) &= \sum_{i\not =j} \phi (\cdot, i, j), \label{gotd:3.5}\\
\lam \eps^2 \phi (\cdot, i, j) +\eps  \phi '(\cdot, i, j) &= \sum_{k\in\mc K} q^j_{k,i}\phi (\cdot, k,j), \quad i\not =j; i\in \mc K.
\label{gotd:4}\end{align} 
These equations imply that $\phi'(\cdot, i,j)$s are absolutely continuous with derivatives in $\alg$. Also, summing, for $j$ fixed, all the equations corresponding to $i\not =j$ we see that 
\( (\lam \eps^2+ 1 ) \sum_{i\not =j} \phi (\cdot, i, j) + \eps \bigl (\sum_{i\not =j} \phi (\cdot, i, j) \bigr )' = \el \phi (\cdot, j, j);\) in this calculation the fact that $Q^j$ is a symmetric matrix satisfying \eqref{sco:1} is used.
Then, inserting the expression for  $\sum_{i\not =j} \phi (\cdot, i, j) $ from \eqref{gotd:3.5} into the so-obtained relation yields, with a bit of algebra, 
\((\ka-1) \phi ''(\cdot, j, j) + \lam \eps (\ka -2) \phi'(\cdot, j, j) - (\lam^2 \eps^2 +\lam \ka )\phi (\cdot, j, j) =0\). The characteristic equation for the so-obtained linear ODE with constant coefficients has two distinct roots, one positive and one negative, the latter being equal to $-\mu$ for $\mu$ defined in \eqref{gotd:1}. It follows that all solutions of the ODE that belong to $\alg$  are of the form given in the first line of \eqref{gotd:3}. 

Next, let $\mathsf 1 \coloneqq(1,1,\ldots,1) \in \R^{\ka-1}$ and  $j$ be still fixed. Eq. \eqref{gotd:4} can be written as 
 an evolution equation (in `time' $x\ge 0$) for the row vector $ (\phi (\cdot, i, j))_{i \neq j}$:
 \begin{align*}
 (\phi (\cdot, i, j))_{i \neq j}' &= (\phi (\cdot, i,j))_{i \neq j} ({\textstyle \frac 1 \eps}\widetilde{Q}^j-\nu {I}) +{\textstyle \frac 1 \eps}\phi (\cdot, j, j)\mathsf 1,\end{align*}
where $I$ is now the $(\ka -1)\times (\ka -1)$ identity matrix, and $\nu$ was defined before the lemma. Hence, by the already established part, for $x\geq 0,$
\begin{align*} 
(\phi (x,i,j))_{i \neq j} & =   (\phi(0,i,j))_{i\neq j}\e^{x(\frac 1\eps\widetilde{Q}^j-\nu {I})} + \tfrac {E_{j,j}}{\eps}\int_0^x   \e^{-\mu y} \mathsf 1  \e^{(x-y)(\frac 1\eps\widetilde{Q}^j-\nu {I})}\ud y. \end{align*}
In other words, introducing $E_{i,j}\coloneqq\phi(0,i,j)$, we have
\begin{align*} 
\phi (x,i,j) & = \e^{-\nu x} \sum_{k\neq j} E_{k,j}   \tilde{p}^j_{k,i}\left(\tfrac{x}{\eps}\right)+\tfrac {E_{j,j}  \e^{-\nu x}}{\eps}\int_0^x \e^{(\nu-\mu) y} \sum_{k\neq j}\tilde{p}^j_{k,i} \left(\tfrac {x-y}{\eps}\right)\ud y, \end{align*}
for $i\not =j$. Since $
 \e^{t \widetilde{Q}^j}$ is doubly stochastic, the sum in the integrand is $1$, and this yields \eqref{gotd:3}. Finally, using the relation between $\mu$ and $\nu$ shown under \eqref{gotd:1} and the fact that 
 $
 \e^{t \widetilde{Q}^j}$ is a stochastic matrix, we check that the functions defined by \eqref{gotd:3} solve \eqref{gotd:3.5}--\eqref{gotd:4} iff \( \sum_{i\not =j} E_{i,j} =\frac{(\ka-1) E_{j,j}}{\eps(\nu-\mu)}, j \in \mc K$.\end{proof}

Our second lemma discusses properties of the approximation \eqref{gotd:7} that constitutes a key to our argument. The following notations are used in the statement of the lemma. For a $\psi \in \Psik$ and a real $\kapp$ matrix $m = \left (m_{i,j}\right )_{i,j\in \mc K}$ we write $\psi m$ to denote $\left (m_{i,j}\psi (\cdot ,j)\right )_{i,j\in \mc K}\in \Phi$. In particular, we will work with the matrices $u,v$ and $w$ defined as follows: 
\begin{align*}
u &= \left (u_{i,j}\right )_{i,j\in \mc K} \, \text{ where } u_{i,j} = 1 \text{ for }  i,j\in \mc K,\\
v &= \left (v_{i,j}\right )_{i,j\in \mc K} \, \text{ where } \, \, v_{i,i} = \el 
\text{ and } v_{i,j}=-1 \text{ for } i\not =j; i,j\in \mc K, \\
w &= \left (w_{i,j}\right )_{i,j\in \mc K} \text{ where } w_{i,i} = \el^2 
\text{ and } w_{i,j}=1 \text{ for } i\not =j; i,j\in \mc K. 
\end{align*}
Also, let $\psi \in \Psik$ be such that $\psi (\cdot, j)$ is absolutely continuous with derivative in $\alg$ for all $j\in \mc K.$ We will write $\psi' $ to denote $ (\psi'(\cdot, j))_{j\in \mc K}$. We note the following relations, which can be easily checked: 
\begin{align} 
\mathfrak A (\psi u) & = \psi' v,\phantom 0 \qquad \mathfrak A (\psi v) = \psi' w, \nonumber \\
\mc Q (\psi u) & = 0,  \phantom{\psi'v} \qquad \mc Q (\psi v) = -\ka \psi v,\nonumber \\
\mc P (\psi v) & = 0,  \phantom{\psi'v} \qquad \mc P (\psi w) = (\ka-1) \psi u,  \label{gotd:6}  
\end{align}
where $\mc P$ was defined in Theorem \ref{thm:2}.

\begin{lemma}\label{lem:gotd2} Let $\psi\in \Psik$ be such that $\psi (\cdot, j)\in W^{3,1}(\R^+), j \in \mc K$, that is, for all $j\in \mc K$,  $\psi (\cdot, j)$   
is three times differentiable with $\psi'(\cdot, j), \psi''(\cdot, j)$ and $\psi ''' (\cdot, j)$ in $\alg$. For $\eps >0$ we define 
\begin{equation}\label{gotd:7} \phi_\eps \coloneqq \psi u + \left ( \eps \kap^{-1} \psi' + \eps^2 \kap^{-2} \psi'' \right )v \in \Phi.\end{equation}  
Then $\phi_\eps \in \dom{\ma} $,  $\grae \phi_\eps = \psi u$  and $\kap \grae \mg \phi_\eps = \psi'' w -  \psi'' v $. Furthermore, $\mc P ( \psi'' w - \psi'' v ) =  (\ka -1) \psi'' u$. 
\end{lemma}
   
\begin{proof} 
The first two claims are immediate. Turning to the third one, we note that, by \eqref{gotd:6}, 
\begin{align*} \eps^{-1} \mathfrak A \phi_\eps &= \eps^{-1} \psi' v + (\ka^{-1} \psi'' + \eps \ka^{-2} \psi''') w \\
\eps^{-2} \mc Q \phi_\eps & = - (\eps^{-1}\psi' + \ka^{-1} \psi'' )v . \end{align*}
It follows that $\ka \grae \mg \phi_\eps = \psi'' w - \psi'' v $, as claimed. The rest is clear by the last two relations in \eqref{gotd:6}. 
\end{proof}

In our third lemma, we explain how transition probability matrices satisfying \eqref{tsos:1} are related to the approximation defined in \eqref{gotd:7}. As a preparation, we consider $F\colon \dom{\ma}\to \R^{\ka (\ka-1)}$ given by 
\begin{equation}\label{tsos:3} \phi \mapsto \left ( \phi (0,i,j) - \el p_{i,j} \phi (0,i,i) -\el  p_{j,j}r_{j,i} \phi (0,j,j) \right)_{j\not =i, i,j\in \mc K},\end{equation}
and note that a $\phi \in \dom{\ma}$ belongs to $\dom{\mc A}$ iff $F\phi =0.$ 

\begin{lemma} \label{lem:tsos1} For $\psi \in \dom{A_\alpha^2}$, let $\phi_\eps , \eps >0$ be defined by \eqref{gotd:7}. Then the limit 
$\grae \eps^{-1} F\phi_\eps $ exists for all such $\psi$ and is finite iff conditions \eqref{tsos:1} are satisfied. In this case 
\[ \grae \eps^{-1} F\phi_\eps =- \textstyle {\frac 1 \kap } F_\alpha \psi \]
where $F_\alpha\colon\dom{A_\alpha}\to \R^{\ka (\ka-1)}$ is given by 
\begin{equation}\label{tsos:fzero} F_\alpha \psi = \big ( D_j(1+\el^2 p_{j,j}r_{j,i}) + D_i \el^2 p_{i,j} \big)_{i\not =j, i,j\in \mc K},\end{equation}
and, see the definition of $\dom{A_\alpha}$, $D_i = \psi'(0,i), i\in \mc K.$  
 \end{lemma}
 
\begin{proof} 
We have
\begin{align*}
\phi_\eps (0,i,i) &= \alpha_i C + (\ka-1) \left (\tfrac \eps \kap D_i + (\tfrac \eps \kap)^2 \psi''(0,i)\right ),\\
\phi_\eps (0,i,j) &= \alpha_j C -  \left (\tfrac \eps \kap D_j +  (\tfrac \eps \kap )^2 \psi''(0,j)\right ), \quad j\not =i; i,j\in \mc K. 
\end{align*}
Hence, the $(i,j)$th coordinate of $F\phi_\eps$ is 
\[ \left ( (1-\el p_{j,j}r_{j,i})\alpha_j - \el p_{i,j} \alpha_i \right ) C - \tfrac \eps \kap \left (D_j (1+\el^2 p_{j,j}r_{j,i}) + D_i\el^2 p_{i,j}  \right ) + o_{i,j} (\eps) ,\]
where $\grae \eps^{-1} o_{i,j} (\eps )=0$. It follows that $\grae \eps^{-1} F\phi_\eps $ exists and is finite if the coefficient next to $C$ is zero in each coordinate, that is, if conditions \eqref{tsos:1} are satisfied.  
\end{proof}

Here is our fourth and final lemma. 

\begin{lemma}\label{lem:tsos2}Let  $\eps,\lam >0$  be fixed. \begin{itemize} 
\item [(a) ] $\ker (\lam - \mg) $ is isomorphic to the subspace $\mathds E $ of $\R^{\kapa}$ formed of $(E_{i,j})_{i,j\in \mc K}$ such that 
\begin{equation}\label{tsos:4}  \sum_{i\not =j} E_{i,j} =\frac{(\ka-1) E_{j,j}}{\eps(\nu-\mu)}, \qquad j\in \mc K. \end{equation}
The isomorphism $\mc I_{\lam,\eps}: \ker (\lam - \mg) \to \mathds E$ identifies the function $\phi \in \ker (\lam - \mg)$ given by 
\eqref{gotd:3} with the matrix of coefficients $E_{i,j}, i,j\in \mc K$ that determines this $\phi.$ 
\item [(b) ]For any $\upsilon\in \R^{\ka(\ka-1)}$ there is precisely one $(E_{i,j})_{i,j\in \mc K}\in \mathds E$ such that 
for the corresponding $\phi= \mc I^{-1}_{\lam, \eps} (E_{i,j})_{i,j\in \mc K}$ we have $F\phi = \upsilon$. In other words, $F$, as restricted to $\ker (\lam - \mg)$ is injective and surjective. 
\item [(c) ] Denoting $(E_{i,j})_{i,j\in \mc K}$ of point \emph{(b)} by $K_{\lam,\eps} \upsilon$ we have 
\( \grae \eps K_{\lam, \eps} \upsilon =  \mu_0^{-1} (\ef \upsilon) \Pi\),
where $\Pi$ was introduced in \eqref{ergo},
\begin{equation}\label{gotd:2}  \mu_0 \coloneqq \sqrt{ {\textstyle\frac{\ka}{\ka-1}} \lam} =\grae \mu (\eps, \lam), \end{equation}
 and the functional $\ef\colon\R^{ \ka(\ka-1)}\to \R$ maps a $\upsilon = (\upsilon_{i,j})_{i\not =j; i,j\in \mc K}$ to the number $\frac1{\ka -1} \sum_{j\in \mc K}\sum_{i\not =j} \upsilon_{i,j}$.
\end{itemize}
 \end{lemma}
\begin{proof} Point (a) is just a restatement of Lemma \ref{lem:gotd1}.

(b) Let $\upsilon= (\upsilon_{i,j})_{i\not =j; i,j\in \mc K}$. Our task is to show that there is precisely one matrix $(E_{i,j})_{i,j\in \mc K} \in \mathds E$ such that, see \eqref{tsos:3},
\begin{equation}\label{tsos:5} E_{i,j}  - \el p_{i,j} E_{i,i} -\el p_{j,j}r_{j,i} E_{j,j} = \upsilon_{i,j}, \qquad i\not =j; i,j\in \mc K, \end{equation}
It is a characteristic feature of this linear system that each $E_{i,j}$ with $i\not =j$ is involved in only one equation there and thus is uniquely determined by $\upsilon_{i,j}$ and `diagonal' elements $E_{i,i}$ and $E_{j,j}$. Moreover, fixing $j\in \mc K$ and summing \eqref{tsos:5} over $i\not =j$ we obtain, by \eqref{tsos:4}, 
\begin{equation}\label{tsos:6}
(\kappa + 1) E_{j,j} - {\textstyle \sum_{i \in \mc K}} p_{i,j} E_{i,i} = \tfrac1{\ka-1}  {\textstyle \sum_{i \not = j}} \upsilon_{i,j}, \qquad j\in \mc K, \end{equation}
where  
\( \kappa = \kappa (\lam, \eps) \coloneqq {\textstyle \frac 1{\eps (\nu - \mu)}} - 1 =\eps \mu +{\textstyle \frac{\eps^2\lam}{\ka-1}}>0,\)
by the relation following \eqref{gotd:1}.
Hence, we are left with showing that the so-obtained reduced system has a unique solution. 

To this end, we note first that the matrix $P$ induces a Markov operator, denoted in what follows by the same letter, $P\colon \R^\kap \to \R^\kap$, given by $P(\xi_j)_{j\in \mc K}= ( \sum_{i\in \mc K} \xi_i p_{i,j})_{j\in \mc K}$; $\R^\kap$ is here seen as an $L^1$-type space, that is, is equipped with the norm $ \| (\xi_j)_{j\in \mc K}\|= \sum_{j\in \mc K} |\xi_j|.$ Hence, the related exponents $P(t)\coloneqq \e^{t(P-I)}, t\ge 0$ (where $I$ is the $\kapp$ identity matrix) are transition matrices of a continuous-time Markov chain whose skeleton is the discrete-time Markov chain described by $P$. In particular, $P(t)$s are Markov operators in $\R^\kap$, and as such they are contractions. It follows that, for any $\rho >0$ and $(\eta_j)_{j\in \mc K}\in \R^\kap$, the resolvent equation for $P-I$: 
\begin{equation}\label{tsos:7} (\rho +1) (\xi_j)_{j\in \mc K} - P(\xi_j)_{j\in \mc K}= (\eta_j)_{j\in \mc K}\end{equation}
has the unique solution 
\[ (\xi_j)_{j\in \mc K} = (\rho +1 - P)^{-1} (\eta_j)_{j\in \mc K} =\left ( \czn \e^{-\rho t} P(t) \ud t \right ) (\eta_j)_{j\in \mc K}. \]
 Since the system \eqref{tsos:6} can be written as 
\[ (\kappa +1 ) \left (E_{j,j} \right )_{j\in \mc K} - P  \left (E_{j,j} \right )_{j\in \mc K} =  \left ( \tfrac1{\ka-1}  {\textstyle \sum_{i \not = j}} \upsilon_{i,j} \right )_{j\in \mc K},\] 
and is thus seen to be a particular case of  \eqref{tsos:7} with $\rho=\kappa$, we obtain
\begin{equation}\label{tsos:8} \left (E_{j,j} \right )_{j\in \mc K} = \left ( \czn \e^{-\kappa t} P(t) \ud t \right )   \left ( \tfrac1{\ka-1}  {\textstyle \sum_{i \not = j}} \upsilon_{i,j} \right )_{j\in \mc K}. \end{equation}


To summarize, the unique solution to \eqref{tsos:5} is described as follows: the diagonal elements $E_{j,j}, j \in \mc K$ 
are determined by \eqref{tsos:8} and the off-diagonal elements are given by
\begin{equation}\label{tsos:9} E_{i,j} =\upsilon_{i,j} +  \el p_{i,j} E_{i,i} +\el p_{j,j}r_{j,i} E_{j,j}, \qquad i\not =j; i,j\in \mc K. \end{equation}

(c) By assumption, \eqref{ergo} holds. It follows that $\lim_{t\to \infty}P(t)=\Pi$ and this in turn implies $\lim_{\rho \to 0} \rho \czn \e^{-\rho t} P(t) \ud t = \Pi$ also.  Moreover, $\grae \kappa =0$ and $\grae \frac \eps \kappa = \mu_0^{-1}$, where $\mu_0$ is defined in \eqref{gotd:2}. Thus, \eqref{tsos:8} shows that 
\begin{equation}\label{tsos:10} \grae \eps (E_{j,j})_{j\in \mc K} = \mu_0^{-1}\Pi   \left (  \tfrac1{\ka-1}  {\textstyle \sum_{i \not = j}} \upsilon_{i,j} \right )_{j\in \mc K}= \mu_0^{-1} (\ef \upsilon ) \alpha_j,\end{equation}
where $\Pi$ is seen as an operator in $\R^\ka$ defined analogously to $P$.
Finally, multiplying \eqref{tsos:9} by $\eps$ and letting $\eps \to 0$ we obtain
\[ \grae \eps E_{i,j} = \mu_0^{-1} (\ef \upsilon) (\el p_{i,j} \alpha_i + \el p_{j,j}r_{j,i}\alpha_j )= \mu_0^{-1} (\ef \upsilon)\alpha_j,  \quad i\not =j; i,j\in \mc K, \]
because \eqref{tsos:1} holds. This combined with \eqref{tsos:10} proves (c). \end{proof}

\begin{corollary} \label{cor:tsos1} For $\psi \in \dom{A_\alpha^2}$, let $\phi_\eps , \eps >0$ be defined by \eqref{gotd:7}. Then, for $K_{\lam,\eps}$ and $\mc I_{\lam, \eps}$ of Lemma \ref{lem:tsos2}, 
\[ \grae K_{\lam, \eps} F\phi_\eps = 0  \mquad{and} \grae \mc I_{\lam,\eps}^{-1} K_{\lam, \eps} F\phi_\eps=0.\]
\end{corollary}

\begin{proof} Lemmas \ref{lem:tsos1} and \ref{lem:tsos2} (c) reveal that \( \grae K_{\lam, \eps} F\phi_\eps = -\frac{\ef F_\alpha \psi}{\ka \mu_0 } \Pi, \)
and therefore for the first part it suffices to check that $\ef F_\alpha \psi =0$. But, since by assumption $\sum_{j\in \mc K} D_j=0$, we have 
\[ \sum_{j\in \mc K} \sum_{i \not =j} (1+\el^2 p_{j,j}r_{j,i} )D_j = \el  \sum_{j\in \mc K} \sum_{i \not =j}p_{j,j}r_{j,i} D_j = \el^2 \sum_{j\in \mc K} p_{j,j}D_j\]
and 
\( \sum_{i \in \mc K} \sum_{j \not =i}\el^2 p_{i,j} D_i = \el^2 \sum_{i\in \mc K} (1-p_{i,i})D_i = - \el^2 \sum_{i\in \mc K} p_{i,i}D_i.\) This shows $\ef F_\alpha \psi =0$.

For the rest it suffices to show that for any $(E_{i,j})_{i,j\in \mc K}\in \mathds E$, the limit $\grae \mc I^{-1}_{\lam,\eps}  (E_{i,j})_{i,j\in \mc K}$ exists, that is, that the functions given by \eqref{gotd:3}, with $E_{i,j}$s fixed, converge as $\eps \to 0$ (in $\alg$). This is indeed the case: the limit function is given by $\phi (x,i,j)=E_{j,j} \e^{-\mu_0 x}, x\ge 0, i,j\in \mc K$. For, we have \eqref{gotd:2}, $\grae \nu =\infty, \grae \eps (\nu - \mu) = 1$ and the norm of the part of $\phi (\cdot, i,j)$ in \eqref{gotd:3} that involves $\e^{-\nu x}$ is bounded by $E\nu^{-1}$ where $E$ is a constant.  \end{proof}

\vspace{-0.75cm}
\section{Proof of the approximation theorem}\label{pota} 

We are finally ready to prove our main theorem. 


Relation \eqref{sco:2} which says that $\lim_{t\to \infty} \e^{t\mc Q} = \mc P$ (even in the operator topology of $\mc L(\Phi)$) 
  allows us to work in the framework of the singular perturbation
  theorem of T. G. Kurtz (\cite{ethier,kurtzper,kurtzapp} or
  \cite[Thm. 42.1]{knigazcup}). Since $\bigcap_{n\ge 1} \dom{(\wdab)^n}$ is a core for $\wdab$ (see Lemma 1.7 p. 53 in \cite{engel}), so is $\dom{(\widetilde \dab)^2}$, and therefore to prove  Theorem \ref{thm:2} we need
  to show  that 
\begin{itemize}
\item [(i) ] for any $\varphi  \in \dom{(\widetilde \dab)^2}$ there are $\varphi_\eps  \in \dom{\mc G_\eps}$ such that 
$
\grae \varphi_\eps = \varphi$, the limit  $\grae \mc G_\eps \varphi_\eps $ exists and $\mc P (\grae \mc G_\eps \varphi_\eps) = 
c\widetilde \dab \varphi$.  
\item [(ii) ]
for any $\phi \in \dom{\mc A}$ we have $
\grae \eps^2 \mc G_\eps \phi = \mc Q \phi$.
\end{itemize}
However, {(ii)} follows {immediately from \eqref{gotd:0}}, and we are left with establishing {(i)}.

So, let $\varphi \in \dom{(\wdab)^2}$. Then, $\psi = J^{-1} \varphi $ belongs to $\dom{(\dab)^2}$ and we can think of $\phi_\eps \in \dom{\mg}$ defined in \eqref{gotd:7}. We know from Lemma \ref{lem:gotd2} that $\grae \phi_\eps = \psi u$ and $\psi u $ is just another notation for $\ka \varphi$. Moreover, $\grae \mg \phi_\eps $ exists and $\mc P (\grae \mg \phi_\eps) = \frac {\ka-1}{ \ka} \psi'' u = (\ka-1) \varphi''.$ 

Hence, we take a $\lam >0$ and define (see Lemma \ref{lem:tsos2}):
\[ \varphi _\eps \coloneqq {\textstyle \frac 1\kap} (\phi_\eps - \mc I_{\lam,\eps}^{-1} K_{\lam, \eps} F\phi_\eps ), \qquad \eps < \eps_0(\lam).  \]
Since $\mc I_{\lam,\eps}^{-1} K_{\lam, \eps} F\phi_\eps$ is a member of $\ker (\lam - \mg)$ such that $F\mc I_{\lam,\eps}^{-1} K_{\lam, \eps} F\phi_\eps= F\phi_\eps $ we see that $F\varphi_\eps =0$, that is, $\varphi_\eps \in \dom{\mc G_\eps}$. Moreover, by Corollary \ref{cor:tsos1}, $\grae \mc I_{\lam,\eps}^{-1} K_{\lam, \eps} F\phi_\eps =0$, and this shows that $\grae \varphi_\eps = \varphi$. Finally, $\mc G_\eps \varphi_\eps = \frac 1\kap (\mg \phi_\eps - \mg \mc I_{\lam,\eps}^{-1} K_{\lam, \eps} F\phi_\eps)  = \frac 1\kap (\mg \phi_\eps - \lam \mc I_{\lam,\eps}^{-1} K_{\lam, \eps} F\phi_\eps) $ because  $\mc I_{\lam,\eps}^{-1} K_{\lam, \eps} F\phi_\eps$ is an eigenvalue of $\mg $ corresponding to $\lam$. It follows that $\grae \mc G_\eps \varphi_\eps = \frac 1\kap \grae \mg \phi_\eps$ and $\mc P (\grae \mc G_\eps \varphi_\eps) = \frac {\ka-1}{ \ka} \varphi''= c\wdab \varphi,$ as desired. 

\vspace{-0.55cm}
\bibliographystyle{plain} 
\bibliography{../../../bibliografia}

\end{document}